\newif
\ifconfver
\confvertrue        %declaring conference version true
\ifconfver
    \documentclass[10pt,journal,final,twoside]{IEEEtran}
\else
    \documentclass[11pt,draftcls,onecolumn]{IEEEtran}
\fi
\usepackage{xspace,empheq,fancybox,amsmath,amssymb,amsfonts,graphicx,epstopdf,epsfig,syntonly,times,amsthm} 
\usepackage{psfrag,color,bm,array,tabularx}
\usepackage{cite,url,footnote,xspace,syntonly,algorithmic}
\usepackage{verbatim,multirow}
\usepackage[linesnumbered,ruled,vlined]{algorithm2e}
\usepackage[T1]{fontenc}
\usepackage{textcomp}
\usepackage{xcolor}
\usepackage{mathtools}
\usepackage{makecell}
\usepackage{bbm, booktabs, makecell, subcaption}
\newcolumntype{M}[1]{>{\centering\arraybackslash}m{#1}}
\newcolumntype{N}{@{}m{0pt}@{}}
\DeclareMathOperator*{\argmax}{\arg\max}
\def\BibTeX{{\rm B\kern-.05em{\sc i\kern-.025em b}\kern-.08em
    T\kern-.1667em\lower.7ex\hbox{E}\kern-.125emX}}

\newtheorem{remark}{\bfseries Remark}
\input{mysymbol.sty}

\newcommand{\KB}{\color{black}{}}
\newcommand{\KBA}{\color{black}{}}
\newcommand{\KBB}{\color{black}{}}
\newcommand{\KBC}{\color{black}{}}
\newcommand{\KBD}{\color{black}{}}

\begin{document}

%%%%%%%%%%%%%%%%%%%%%%%%%%%%%%%%%%%%%%%%%%%%%%%%%%%%%%%%%%%%%%%%%%%%%%
%                                                                    %
%               Paper Title                                         %
%                                                                    %
%%%%%%%%%%%%%%%%%%%%%%%%%%%%%%%%%%%%%%%%%%%%%%%%%%%%%%%%%%%%%%%%%%%%%%

\title{Reinforcement Learning Based Optimal Battery Control Under Cycle-based Degradation Cost}
\author{
\IEEEauthorblockN{Kyung-bin Kwon},~\IEEEmembership{Student Member, IEEE}, and \IEEEauthorblockN{Hao Zhu},~\IEEEmembership{Senior Member, IEEE}

\thanks{\protect\rule{0pt}{3mm} Manuscript received August 21, 2021; revised January 16, 2022 and April 24, 2022; and accepted May 28, 2022. This work has been supported by NSF Grants 1802319, 1952193, and 2130706.}
%This work has been partially supported by the Illinois Center for a Smarter Electric Grid (ICSEG), Trustworthy Cyber Infrastructure for the Power Grid (TCIPG), and Cyber Resilient Energy Delivery Consortium (CREDC).} 
\thanks{\protect\rule{0pt}{3mm}  The authors are with the Department of Electrical \& Computer Engineering, The University of Texas at Austin, 2501 Speedway, Austin, TX, 78712, USA; Emails: {\{kwon8908kr, haozhu\}{@}utexas.edu}.
}}

% % The paper headers
\markboth{IEEE TRANSACTIONS ON SMART GRID (ACCEPTED)}%
{Kwon \MakeLowercase{\textit{et al.}}: Reinforcement Learning Based Optimal Battery Control Under Cycle-based Degradation Cost}
\renewcommand{\thepage}{}
\maketitle
\pagenumbering{arabic}

%%%%%%%%%%%%%%%%%%%%%%%%%%%%%%%%%%%%%%%%%%%%%%%%%%%%%%%%%%%%%%%%%%%%%%
%                                                                    %
%                   Abstract                                         %
%                                                                    %
%%%%%%%%%%%%%%%%%%%%%%%%%%%%%%%%%%%%%%%%%%%%%%%%%%%%%%%%%%%%%%%%%%%%%%
%
\begin{abstract}
Battery energy storage systems are providing increasing level of benefits to power grid operations by decreasing the resource uncertainty and supporting frequency regulation. Thus, it is crucial to obtain the optimal policy for utility-level battery to efficiently provide these grid-services while accounting for its degradation cost. To solve the optimal battery control (OBC) problem using the powerful reinforcement learning (RL) algorithms, this paper aims to develop a new representation of the cycle-based battery degradation model according to the rainflow algorithm. As the latter depends on the full trajectory, existing work has to rely on linearized approximation for converting it into instantaneous terms for the Markov Decision Process (MDP) based formulation. We propose a new MDP form by introducing additional state variables to keep track of past switching points for determining the cycle depth. The proposed degradation model allows to adopt the powerful deep Q-Network (DQN) based RL algorithm to efficiently search for the OBC policy. Numerical tests using real market data  have demonstrated the performance improvements of the proposed cycle-based degradation model in enhancing the battery operations  while mitigating its degradation, as compared to earlier work using the linearized approximation.
\end{abstract}


%\newpage

%%%%%%%%%%%%%%%%%%%%%%%%%%%%%%%%%%%%%%%%%%%%%%%%%%%%%%%%%%%%%%%%%%%%%%
% %%
% %%      Section: Introduction %%
% %%%%%%%%%%%%%%%%%%%%%%%%%%%%%%%%%%%%%%%%%%%%%%%%%%%%%%%%%%%%%%%%%%%%
%%%

\section{Introduction}\label{sec:intro}

Battery energy storage systems as flexible resources are a key technology to enable the decarbonization of electricity infrastructure  in future \cite{Arbabzadeh2019, denholm2010role}. {\KB Particularly, utility-level battery systems can be used to increase the payoff from electricity market via energy arbitrage \cite{1-3}, while contributing to the grid's power balance through participating in ancillary services \cite{1-4}.}
	% and have the additional profit for providing the services in the unitlity level.} %The owner can not only use the battery to the energy arbitrage to maximize the profit \cite{1-3} but also participate in ancillary services such as voltage regulation, peak-shaving and frequency regulation and thus contribute to the stability of the power system \cite{1-4}. 
It is crucial to develop effective strategies for real-time battery operations in order to utilize its flexibility potentials to mitigate  the increasing uncertainty introduced by renewable or non-controllable loads. 

{\KB The optimal battery control (OBC) problem for determining the (dis)charging policies has been popularly considered to reduce a combination of battery operational costs. It aims to reduce the net cost for electricity usage and frequency regulation (FR) penalty, as well as possible violations of network constraints; see e.g., \cite{UW,2-3,2-4,2-7}. } %Additional problem consideration includes the precise circuit model of battery \cite{2-5} or the network constraints on grid voltage \cite{2-7}.
In addition, battery's cycle life as characterized by the degradation cost is especially needed when participating in FR or other fast services \cite{UW,2-4}. Unlike other costs that mostly depend on the instantaneous battery status, the modeling of battery degradation is \textit{cycle-based} according to the full trajectory of battery's state of charge (SoC). It requires the identification of all charging/discharging cycles using the so-termed \textit{rainflow} algorithm  \cite{4-4}. Thus, degradation-aware OBC problem results in increased complexity as shown by \cite{UW,2-4}.

 {\KB Due to the fast dynamics in prices or load demands,  the OBC solution can be greatly affected by the uncertainty of future information. To address this issue, a model-predictive control (MPC) framework has been widely used by optimizing the current action according to the predicted input values for a fixed time window; see e.g., \cite{morstyn2017model,MPC_LFC,MPC01}. } {\KB Nonetheless, the FR signal exhibits very minimal temporal correlation \cite{FR_ref}, leading to significant difficulty in predicting  it and thus applying MPC for reducing FR penalty. Furthermore, even though battery health has been considered in MPC-based OBC work \cite{MPC1, MPC2}, the cycle-based degradation model is largely missing.} %Moreover, because of the randomness that frequency regulation signal has \cite{FR_ref}, predicting the future information can be a difficult task for MPC.}

Recent advances on reinforcement learning (RL) \cite{RL} enable the effective search of optimal control policies directly using real data samples to address the uncertainty issue in dynamical systems.  This data-driven framework helps to bypass the hurdles in formulating the complex models of system dynamics or estimating the statistical information on the uncertainty. Specifically for the OBC problem, it allows to flexibly incorporate a variety of operational objectives, and several RL techniques have been widely used,  such as the Deep Q-Network (DQN) \cite{3-2}, SARSA \cite{3-3} and $TD(\lambda)$-learning \cite{3-4}.  However, a majority of these techniques have not considered the battery degradation cost, due to the difficulty of representing cycle-based model in the RL formulation. Very recently,  \cite{4-6} has developed a linearized approximation for the cycle-based degradation cost, which converts it to an instantaneous degradation coefficient  that can be used by the RL algorithms. Nevertheless, the accuracy of this approximation method depends on a given sample trajectory based on which the linearization is performed. %{\KB Similarly in \cite{QMPC}, which implements Q-learning to the MPC approach, the battery degradation cost is considered as a linear funnction of the batttery (dis)charging power with the pre-defined coefficient.}
 The resultant modeling mismatch can limit the RL iterations from finding the best policy within the full  search space. Thus, it is still an open problem of effectively incorporating the accurate battery degradation cost into the search of {\KB OBC} policy. 

%For example in \cite{3-1} and \cite{3-2}, the basic MDP formulation of the battery management problem has been made and the reinforcement learning is applied; the battery operation is optimized to maximize the profit by arbitrage in the electricity market using Policy Function Approximation method and Deep Q-Network method, respectively. 
%In \cite{3-3}, the control of the battery operation includes buying power, meet the demand and provide frequency regulation service and the Sarsa algorithm is applied to optimize the long-term rewards based on meeting the load demand. In \cite{3-4}, standard batch RL technique is applied to minimize the battery operation cost. Both the battery's charge and discharge efficiencies and the nonlinearlity model of the inverter's efficiency in the microgrid are considered when generating MDP model.  Similar to \cite{3-4}, the residential energy storage system contol is proposed using RL for electric bill minimization. $TD(\lambda)$-learning algorithm is applied when minimizing the battery charging cost, which provides the higher convergence rate and higher performance. In \cite{3-6}, the energy manamement algorithm is proposed for smart energy buildings considering both the battery operation and V2G (Vehicle-to-Grid) system of EVs. The common features of the papers are that the uncertainty from the power system has been included such as demand, electricity price and renewable generation.

The goal of our work is to develop a modeling approach to precisely represent the battery degradation cost and use it for the design of RL-based OBC algorithm. {\KB The overall objective includes the net electricity cost, FR penalty, and cycle-based degradation cost.} The main modeling challenge lies in the latter as it is determined by the battery's full SoC trajectory. Based on the rainflow algorithm, the complex process of material fatigue is associated with the stress level of each individual charging or discharging cycle \cite{4-5}. Thus, the degradation cost is an exponentially increasing function of cycle depth \cite{4-4}, and the latter strongly depends on the past trajectory of battery status.  This leads to a pronounced mismatch with the Markov Decision Process (MDP) form used by RL algorithms, as the latter would represent the problem objective as functions of instantaneous states and actions only. The aforementioned approach of linearizing the degradation cost as in \cite{4-6,4-3} fails to recognize this exponential relation with the cycle depth, and unfortunately can lead to deep (dis)charging cycles that may not be overall profitable. 

To this end, we have analytically shown that it is possible to keep track of the battery cycles by augmenting the state with the more recent switching points (SPs) along the SoC trajectory. These critical transition points between charging and discharging sessions are extremely useful for identifying the correct cycle depth according to the rainflow condition. In addition, they allow for decomposing the degradation cost of a full (dis)charging cycle into incremental differences between consecutive time instances in the form of instantaneous cost.  This proposed representation of battery degradation cost helps to deploy state-of-the-art RL algorithms to learn the OBC policy. We have used the DQN technique to search for the parameters of the action-value function, or Q-function, associated with the resultant MDP form. 
 
 To sum up, the main contribution of the present paper is two-fold. First, we have developed an instantaneous cost modeling of the battery degradation with \textit{guaranteed equivalence} to the original cycle-based representation based on rainflow algorithm. Second, the proposed degradation cost is successfully applied to form an MDP, allowing to develop efficient RL algorithms for the OBC problem. Our numerical tests using real data of electricity prices and FR signals  have validated the performance improvement of the proposed cycle-based cost model in accurately representing battery degradation and effectively generating profit-seeking OBC policies.

The rest of the paper is organized as follows. Section \ref{sec:SM} introduces the key variables for modeling the battery control problem into the MDP form. In Section \ref{sec:cost}, we model the cycle-based degradation cost using the rainflow algorithm, and develop a new approach to represent it as instantaneous cost through state augmentation. Section \ref{sec:bc} formalizes the OBC problem and presents the DQN method as the RL solution technique. Numerical results using real-world data are presented in Section \ref{sec:CS} to validate the performance improvement of the proposed degradation model, as compared to earlier approach using linearized approximation. Finally, the paper is wrapped up in Section \ref{sec:CON}. %{\KBC with one appendix.}

\begin{table}[t]
    \centering
        \caption{List of symbols}
    \label{tb_symbols}
\begin{tabular}{cl}
\toprule
\textbf{Notation} & \multicolumn{1}{c}{\textbf{Description}} \vspace*{3pt} \\ \hline \hline
$c_t$ & state of charge (SoC) of the battery \\
$\overline{c}, \underline{c}$ & maximum/minimum capacity of the battery \\
$p_t$ & electricity market price \\
$f_t$ & frequency regulation (FR) signal \\
$s_t$ & battery full state \\
$a_t$ & battery charging/discharging action\\
$b_t$ & battery charging/discharging power \\
$\gamma$& discount factor\\
$\ccalT$ & the exploration time-horizon\\
$h^e_t$ & net cost for electricity usage \\
$h^f_t$ & frequency regulation penalty \\
$h^d_t$ & battery degradation cost \\
$\delta$ & frequency regulation penalty coefficient \\
$\Phi(d)$ & degradation cost for a cycle of depth $d$ \\
$\alpha_d, \beta$ & degradation coefficient based on battery types \\
$c^{(0)}_t, c^{(1)}_t, c^{(2)}_t$ & SoC level of switching points (SPs) \\ \bottomrule
\end{tabular}
\end{table}

%%%%%%%%%%%%%%%%%%%%%%%%%%%%%%%%%%%%%%%%%%%%%%%%%%%%%%%%%%%%%%%%%%%%%%
% %%
% %%      Section: Problem Formulation %%
% %%%%%%%%%%%%%%%%%%%%%%%%%%%%%%%%%%%%%%%%%%%%%%%%%%%%%%%%%%%%%%%%%%%%
%%%

%%%%%%%%%%%%%%%%%%%%%%%%%% section 1 %%%%%%%%%%%%%%%%
\section{System Modeling} \label{sec:SM}

{\KB This paper considers the optimal battery control (OBC) for maximizing the economic pay-off while accounting for the battery degradation. The pay-off is from energy market participation and also the provision of FR service, as discussed later.  One notable feature of the present work is the consideration of battery degradation cost, which can greatly increase the life-cycle under any  general pay-off model {\cite{UW}}. A list of symbol notation and description is tabulated in Table \ref{tb_symbols}.}

To determine the battery's {effective (dis)charging power} $b_t \in [\underline{b},~\bbarb]$ at each discrete-time instance $t = 0, 1, \ldots$, we introduce a list of state variables based on battery status or external inputs. 
\begin{itemize}
	\itemsep 2pt
	\item $c_t \in [\underline{c},~\bbarc]$: normalized state of charge (SoC) of the battery;
\item $p_t$: electricity market price; 
\item $f_t$: frequency regulation (FR) signal. 
\end{itemize}
Note that the SoC is normalized by the maximum capacity; i.e., $c_t \in [0,1]$. It is also an internal battery state affected by the past actions $\{b_\tau \}$, whereas the other states are received from grid operators and thus are not directly action-dependent. 

To leverage reinforcement learning algorithms for this problem, we consider it as a Markov Decision Process (MDP) {\cite[Ch.~3]{MDP}} denoted by a tuple $(\mathcal{S}, \mathcal{A}, \mathcal{P}, \mathcal{R}, \gamma)$, as detailed here. % Specifically, $\mathcal{S}$ and $\mathcal{A}$ are the state and action spaces, respectively. The transition kernel  $\mathcal{P}: \mathcal{S} \times \mathcal{A} \times \mathcal{S} \rightarrow [0,1] $ refers to the system dynamics, while the reward function $\mathcal{R}: \mathcal{S} \times \mathcal{A} \rightarrow \mathbb{R}$ captures the learning objective and $\gamma \in (0,1)$ is a discount factor, as detailed here.

\textbf{State space} $\mathcal{S}$ contains the set of feasible values for the system state $s_t$, including both the SoC $c_t$, and the other inputs $p_t$ and $f_t$ which affect the economic benefits. {\KB Additional state variables will be specified in Section~\ref{sec:bc} for representing cycle-based degradation cost.} State dynamics need to follow the Markov property as discussed soon. 

\textbf{Action space} $\mathcal{A}$ includes the set of decisions that battery can take. {\KB We consider a discrete multi-level set with a total of $|\ccalA|$ actions, as 
 \begin{align}
 a_t \in \mathcal{A}=\{a^{(1)}, a^{(2)}, \cdots, a^{(|\ccalA|)} \}  \label{action2}
 \end{align}
with normalized actions $a^{(n)} \in [-1,~1]$.  Accordingly, the normalized (dis)charging power $b_t \in [ \underline{b},~\bbarb]$ is set to be
 \begin{align}
 b_t = \begin{cases}
 \min\{\overline{c}-c_t, \overline{b} a_t\}  &\;\text{if}\; a_t \geq 0,\\ %\overline{b} a_t \geq \overline{c}-c_t\\
% \overline{b} a_t &\;\text{if}\; a_t  \geq 0,  \overline{b} a_t \leq \overline{c}-c_t\\
\max \{  \underline{c}-c_t,  \underline{b} a_t\} &\;\text{if}\; a_t < 0.
%, \underline{c}-c_t \leq \underline{b} a_t\\ &\;\text{otherwise}
 \end{cases}  \label{bt}
 \end{align}}
Continuous action space that directly determines $b_t=a_t$ is also possible. While this paper focuses on a discrete $\ccalA$, the RL algorithm can be generalized to continuous $a_t$ as well. 

\textbf{Transition kernel} $\mathcal{P}: \mathcal{S} \times \mathcal{A} \times \mathcal{S} \rightarrow [0,1] $ captures the system dynamics under the Markov property {\cite[Ch.~3]{MDP}}. For the input states such as price $p_t$, we assume Pr$(p_{t+1}|\{p_\tau\}_{\tau=1}^t)=$ Pr$(p_{t+1}|p_t)$; and similarly for $f_t$. This is reasonable as the market price has very short-term memory {\cite{MPM}}, while FR signal $f_t$ can be modeled as a white noise sequence of no memory {\cite{FR}}. A longer memory is possible too; such as the prices that follow Pr$(p_{t+1}|\{p_\tau\}_{\tau=1}^t)=$ Pr$(p_{t+1}|p_t,~p_{t-1})$. In this case,  both  $p_t$ and $p_{t-1}$ are included as the part of the state per time $t$ to satisfy the Markov transition property.
%{\hao why?? \cite{PriceM}}{\KB (reference changed)}.   

Using Eq.~\eqref{bt}, the  SoC state $c_t$ transitions as
\begin{align}
 c_{t+1} = c_t+b_t,~\text{with}~b_t~\text{given in}~\eqref{bt}.
\label{transition_c}
\end{align}
For general action space with any  $b_t\in[\underline{b},\bbarb]$, $c_t$ is updated by
\begin{align}
c_{t+1} = \begin{cases}
\overline{c} \; &\text{if}\; \overline{c}-c_t \leq b_t,\\
\underline{c} \; &\text{if}\; \underline{c}-c_t \geq b_t,\\
c_t+b_t \; &\text{otherwise}.
\end{cases} \label{transition_cc}
\end{align}

\textbf{Reward} function $\mathcal{R}: \mathcal{S} \times \mathcal{A} \rightarrow \mathbb{R}$ captures the learning objective. Notably, it is always the accumulated reward consisting of \textit{instantaneous} terms, where per time $t$  the latter only depends on the current state and action as
 \begin{align}
 %\mathcal{R}: 
 r_{t}=r_t(s_t, a_t) \label{reward}
 \end{align} 
In the following we will minimize the objective cost function $h_t$ as negative reward, where its instantaneous property will be ensured after introducing additional state variables {\KB as detailed in Section~\ref{sec:bc}.}
%When computing and keeping the instantaneous reward , more state variables may need to be considered according to the cost functions. The detailed description of cost would be discussed in Section \ref{subsec:cost}.\\
 
\textbf{Discount factor} $\gamma \in (0,1]$ is a constant to accumulate the total reward along the time horizon. Smaller $\gamma$ values imply that future rewards are less important than current ones at a discounted rate \cite[Ch.~3]{MDP}. As we adopt a finite exploration time-horizon $\mathcal T = [1,\ldots,T]$  for the OBC problem, for simplicity $\gamma = 1$ will be used. 

%%%%%%%%%%%%%%%%%%%%%%%%%% section 3 %%%%%%%%%%%%%%%%
\section{Modeling of Battery Degradation Cost} \label{sec:cost}

We consider three types of operational cost related to battery management. The energy cost relates to the electricity price according to (dis)charging, while the FR cost is based on its fast-varying flexibility. Under a contract of providing FR service, the battery would follow the $f_t$ signal sent by the market operator as much as possible {\cite{UW}}. These two costs can be simply obtained by the state variables discussed so far. First, the net cost for electricity usage under (dis)charging power $b_t$ can be represented as
\begin{align}
h^e_t(p_t, b_t) = p_t b_t, \quad \forall~t \in \ccalT. \label{cost_e}
\end{align}
Second, using a penalty coefficient $\delta$ for deviation from FR signal $f_t$, one can form 
\begin{align}
h^f_t (f_t, b_t) = \delta \lvert f_t - b_t \rvert,  \quad \forall t \in \ccalT. \label{cost_f}
\end{align}
{\KB The energy cost in \eqref{cost_e} is typically negative due to the energy arbitrage capability, while the FR penalty in \eqref{cost_f} is always positive. This is because the additional economic benefit by participating in the FR contract is not included here. Overall, a battery should receive positive pay-off from these two tasks. }

\begin{remark} \label{rmk:frsignal}\emph{(Frequency regulation signal)}
	{\KB In practice, the FR signal is much faster than other system dynamics. For example, the real-time price is typically updated every 5 minutes, while the FR signal may be at 2-second rate \cite{pjm_fr}. To reduce the complexity of the training computation later on, we will down-sample the FR signal to attain $\{f_t\}$ at a slower rate for searching the policy. In testing and implementing the resultant  policy, the original fast FR signal will  be instead used to realistically evaluate the performance of the RL approach.}
	\end{remark}

{\KBB As for the battery degradation cost, there are several stress factors affecting the battery lifetime such as temperature, {\KBC high C-rates}, average SoC, and Depth of Discharge (DoD) \cite{4-4,wang}. During daily battery operations, the DoD stress model is considered the most relevant while other factors may be minimally affected.  {\KBC This will be shown numerically in Section~\ref{sec:CS}}.  
	%Among different stress factor models, we focus on implementing  stress model which has close relation to daily real-time operation of the battery. 
According to the DoD stress model, the aging of battery cells mainly depends on material fatigue as a result of (dis)charging cycles of the SoC trajectory, especially due to following the FR signal {\cite{FR_ref}}. Since this cycle-based degradation constitutes as a key battery lifetime consideration \cite{forementioned}, the proposed OBC formulation to reduce it can greatly increase the battery's lifetime revenue.}

%Compared to the two cost functions previously expressed using instantaneous state/action variables, the \textit{cycle-based} degradation cost is more complicated to model. The key challenge becomes to effectively decompose it into instantaneous terms for each $t$. 
The rainflow algorithm \cite{4-4} is widely used for computing the cycle-based degradation cost. 
%{\hao (need a paper from battery research on rainflow, not another battery control paper that used rainflow)}{\KB (changed to the paper related to battery research on rainflow)} 
%{\hao (start of new edits)} 
\begin{figure}[t]
	\centering
	\includegraphics[width=\linewidth]{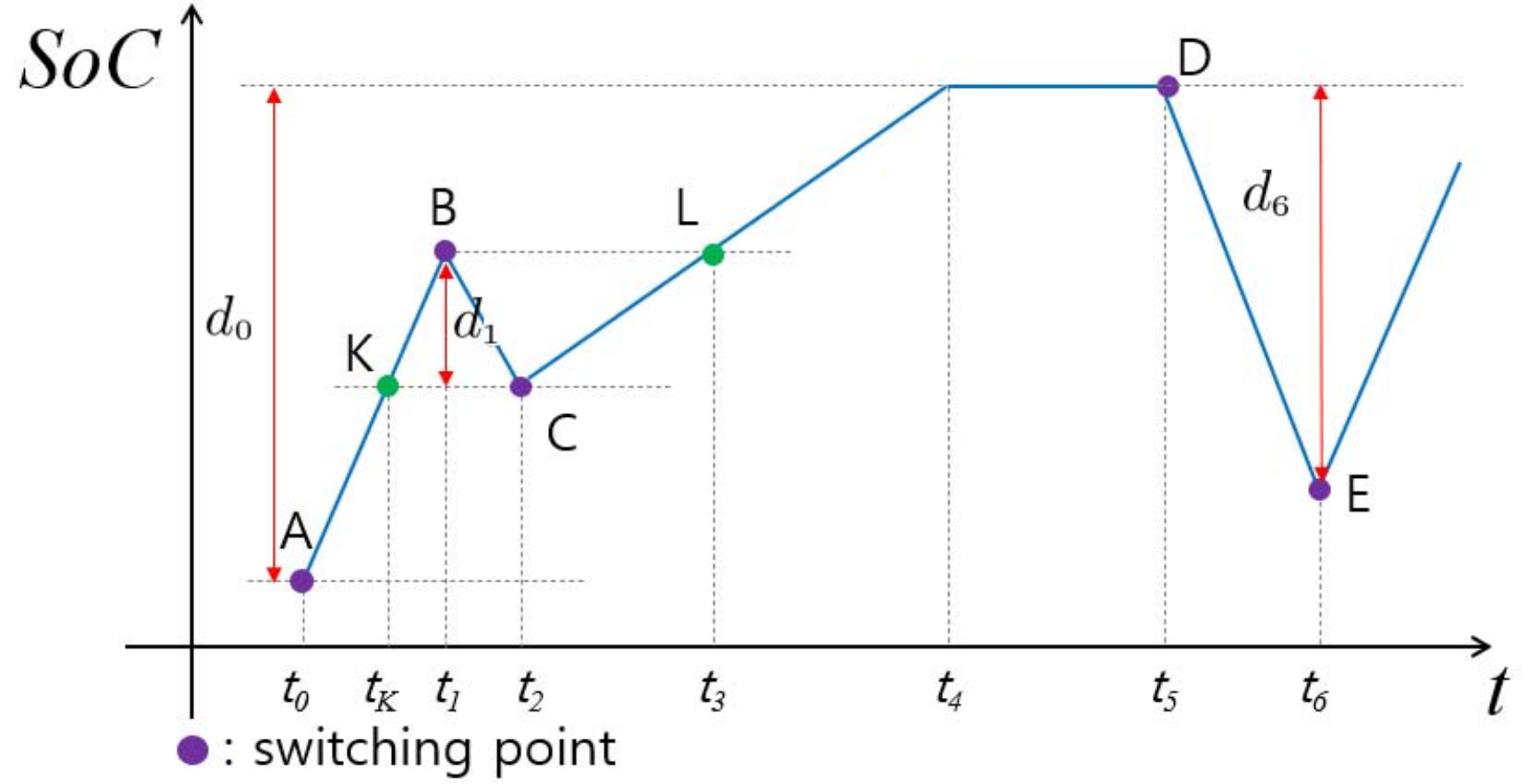}
	\caption{An example of battery SoC trajectory used for modeling the battery degradation cost based on the rainflow algorithm.} %{\hao (why do we need t in the symbols $d_t0$? can we just have $d_0$?? )}{\KB (changed to have only the number)}}
\label{fig_rainflow_example}
\end{figure}
Fig.~\ref{fig_rainflow_example} illustrates an example of battery SoC trajectory which consists of several charging and discharging cycles. The  switching points (SPs), labeled by $A-E$, correspond to the transitions between charging and discharging and will be used for identifying the cycles by rainflow algorithm. %Rainflow algorithm is a method to identify the cycle in material fatigue analysis, including the battery degradation \cite{rainflow2}, \cite{rainflow3}. 
For example, the trajectory $A-B-C-D$ consists of a long charging cycle with a small discharging part from $B-C$. The respective depths of these two cycles, defined as the absolute SoC differences between the start and end SPs, are $d_0$  and $d_1$. As $d_1$ is smaller than the difference between $A-B$ and that between $C-D$, this trajectory is thus divided into the \textit{full cycle} from $K-B-C$ of depth $d_1$ and the other \textit{half cycle} from $A-K(C)-D$ of depth $d_0$. This is the so-called rainflow condition as stated in Lemma \ref{lemma:sp}; {\KBC see e.g., \cite{UW}.}

%as formally stated here.
% which uses \textit{four consecutive SPs} to determine the exact cycles and their depths, .

{\KBB
\begin{lemma} \label{lemma:sp}
		The SoC values of the last three SPs by time $t$ are sufficient for evaluating the rainflow condition and determining the depth of (dis)charging cycles.
		% based on the rainflow condition by comparing the SoC level differences between the neighboring SPs and current SoC.
\end{lemma}
}

{Based on the cycle depth $d>0$, the associated degradation cost is given by} 
\begin{align}
\Phi (d) = \alpha_d e^{\beta d} \label{cycle_cost}
\end{align}
with positive constant coefficients $\alpha_d$ and $\beta$ based on battery types{\KBC \cite{4-3,4-4}}. Recalling the normalized SoC $c_t \in [0,1]$, we have the cycle depth $d\in[0,1]$ as well. {Note that for any pair in $\ccalD:= \{(d_1,d_2):~d_1,~d_2 \geq 0,~d_1+d_2 \leq 1\}$, we can show that $e^{(d_1+d_2)} \leq e^{d_1} + e^{d_2}$. This is because the maximum value of the function $g(d_1, d_2) := e^{(d_1+d_2)}-e^{d_1}-e^{d_2}$ for the simplex $\ccalD$ equals to $(e-2e^{0.5}) <0$, which is attained at $(d_1, d_2) = (0.5, 0.5)$. Thus, to reduce the degradation cost a single (dis)charging cycle that is longer and deeper is typically preferred, as opposed to the combination of multiple shorter cycles. This intuitive rule for cycle-based degradation model will be demonstrated later on in numerical tests. } Unfortunately, this cycle-based degradation cost depends on  the past SoC trajectory, and unfortunately, it does not follow the accumulated form of instantaneous terms as in Eq.~\eqref{reward}. %To address this issue, existing approach approximates it using a linear model, as detail here.

{\KB \textbf{Linearized degradation model} has been developed in \cite{4-6}  to compute the averaged degradation coefficient from past SoC trajectory. 
	%and approximates the instantaneous cost caused by each (dis)charging action.
Specifically, a degradation coefficient $\alpha_d$ is first determined using a given SoC trajectory over $\ccalT$ as} %($\bar{\cdot}$ indicating sample)
\begin{align}
	a_{d} = \frac{\sum_{i=0}^{\bar N} \Phi (\bbard_i)}{\sum_{t=0}^T \lvert \bbarb_t \rvert} \label{deg_coeff}
\end{align}
%
%{\hao (since we haven't introduced the $h^d$ why having these terms show up??? what is the subscript j mean??)}
by averaging the total degradation costs of the $(\bbarN+1)$ cycles over the accumulative absolute charging power throughout the sample trajectory. This way, the instantaneous degradation cost for any new SoC trajectory is approximated by  
\begin{align}
h^{d}_t (b_t) \approxeq -a_d \lvert b_t \rvert. \label{linear_deg}
\end{align}
This linearized degradation cost model can be easily computed once $a_d$ is known. However, this approximation inexplicitly assumes that the new trajectory should be very similar to the given sample trajectory for computing $a_d$. To implement the RL algorithm later on, the coefficient $a_d$ will be updated using the most recent trajectory during the sampling process. Nonetheless, as an approximation it does not represent the actual cycle-based degradation cost and thus limits the RL algorithm's search for the best SoC trajectory. 

\begin{figure}[t]
	\centering
	\includegraphics[width=.85\linewidth]{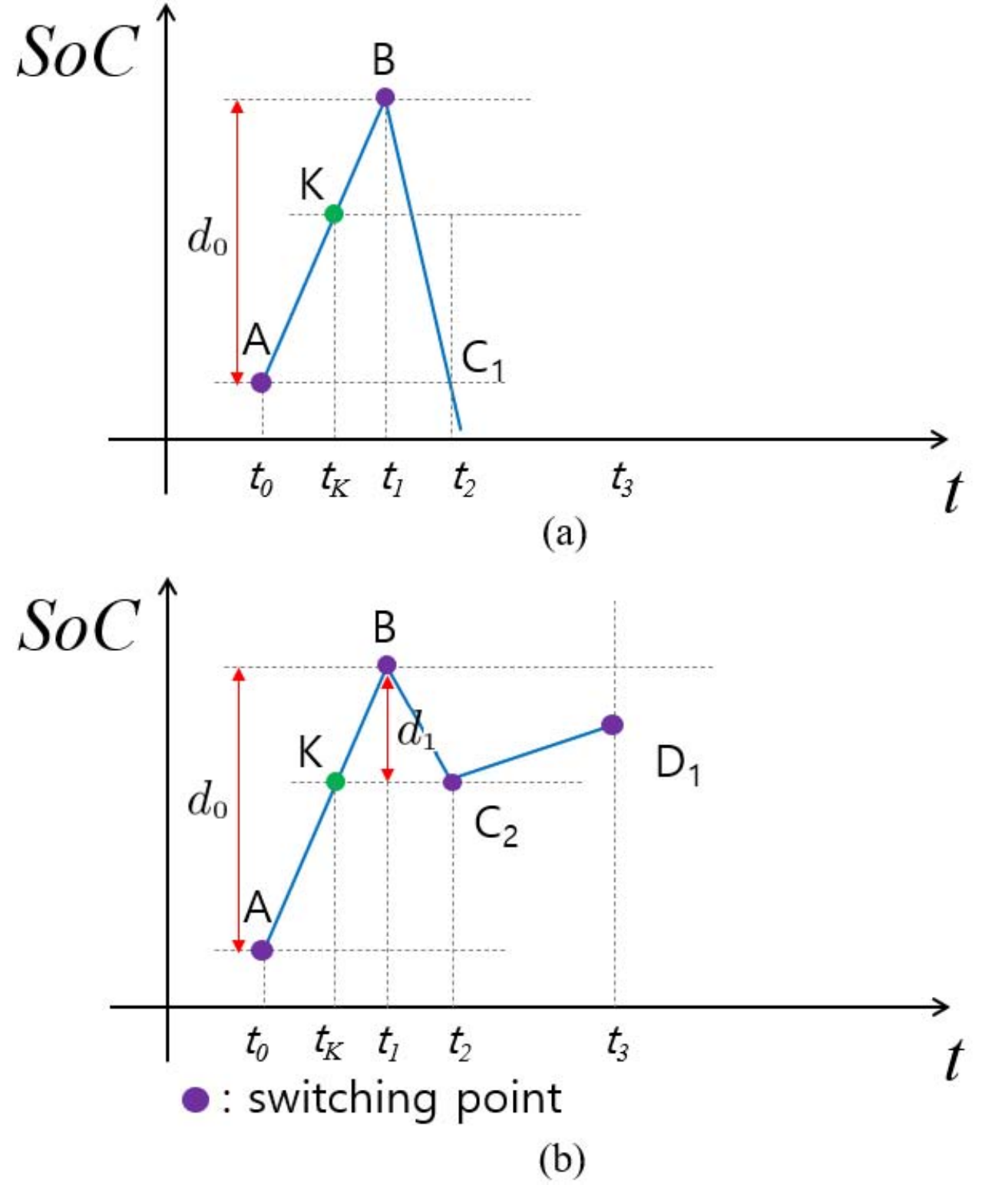}
	\caption{Two cases of rainflow condition not satisfied: (a) case $NR_a$ and (b) case $NR_b$.} %{\hao (a lot of empty space at the top; why the first depth is $d_1$ not $d_0$??)}{\KB (changed the number)}}
	\label{fig_no_rainflow_cases}
\end{figure}

\textbf{Cycle-based degradation model} will be pursued instead to address the approximation issue by augmenting the state $s_t$ with the last three SPs before time $t$. As stated in Lemma \ref{lemma:sp}, they are sufficient information for checking the rainflow condition. The state $s_t$ now includes three additional variables, $c^{(0)}_t$, $c^{(1)}_t$, and $c^{(2)}_t$, as the SoC from the oldest SP to the latest one. Note that they may overlap if there are less than three SPs before time $t$. For example, at point $K$ in Fig.~\ref{fig_rainflow_example},  these three SP states all equal to the SoC of point $A$; and similarly for point $C$, we have $c^{(1)}_t=c^{(2)}_t$ equal to the SoC of $B$. The latest SP's SoC $c^{(2)}_t$ can be used to identify if the current instance $t$ is a new SP, using the rule %by comparing with the sign of $b_t$, as
\begin{align}
b_t (c_t-c^{(2)}_t) < 0. \label{rainflow_point}
\end{align}
If Eq.~\eqref{rainflow_point} holds, we have a new SP and will update $\{c^{(i)}_t\}$ based on whether the rainflow condition is satisfied. 

%In addition, the SPs are updated according to whether the rainflow condition is satisfied at time $t$. 
To update $\{c^{(i)}_t\}$, Fig.~\ref{fig_no_rainflow_cases} illustrates two cases where the rainflow condition is not satisfied.  Fig.~\ref{fig_no_rainflow_cases}(a) shows the SoC of the point $C_1$ is not within the range between $A$ and $B$, while Fig.~\ref{fig_no_rainflow_cases}(b) indicates the SoC of the new SP $D$ is within the range between $C_2$ and $B$. These cases are denoted by cases $NR_a$ and $NR_b$, respectively. In either case, {\KBB the oldest SP $A$  will be removed while the remaining SPs will be used to update $\{c^{(i)}_{t+1}\}$, as listed in Table \ref{tb_transition}.} In addition, case $RA$ denotes the scenario of rainflow condition being satisfied, such as the point $D$ in  Fig.~\ref{fig_rainflow_example}.  This is because at SP $D$: i) the third SP $C$ lies between the first two SPs $A$ and $B$; and ii) the current SoC at SP $D$ exceeds the range between the latest two SPs $B$ and $C$. Hence, the trajectory $A-B-C-D$ is divided in to the long half-cycle $A-K-C-L$, and another full-cycle $K-B-C$ of depth $d_1$. After the $RA$ case is satisfied, the two SPs $B$ and $C$ will be removed from the record. The SoC state updates for all three cases are summarized in Table \ref{tb_transition}.

Interestingly, the state transitions in Table \ref{tb_transition} also allow for decomposing the cycle-based degradation cost into instantaneous difference term for each instance $t$. As the cycle depth changes according to $b_t$ only, the degradation cost in Eq.~\eqref{cycle_cost} can be modeled by accumulating the following incremental term per time $t$:
\begin{align}
h^{d}_t(b_t, c_t, c^{(2)}_t) = \alpha_d e^{\beta \lvert c_t+b_t-c^{(2)}_t\rvert}-\alpha_d e^{\beta \lvert c_t-c^{(2)}_t \rvert}. \label{inst_deg_cost}
\end{align}
Basically, the instantaneous degradation model in Eq.~\eqref{inst_deg_cost} accounts for difference of $\Phi(\cdot)$ due to the change of cycle-depth, which can be computed based on the latest SP state $c^{(2)}_t$. {\KBC Therefore, summing up all instantaneous terms in Eq.~\eqref{inst_deg_cost} yields the total degradation cost, as formally stated in \textbf{Proposition 1} with the proof provided in the Appendix.} 
{\KBA
\begin{proposition} \label{prop:action}
 Under Lemma \ref{lemma:sp},  the summation of the instantaneous terms in Eq.~\eqref{inst_deg_cost} throughout the time-horizon $\mathcal{T} = [1, \cdots, T]$ is exactly equivalent to the total cycle-based degradation cost along $\mathcal{T}$. 
\end{proposition}
}

\begin{table}[t]
    \centering
    \caption{State transitions at a new SP identified at time $t$}
    \label{tb_transition}
\begin{tabular}{c|ccc}
        \Xhline{2\arrayrulewidth}
\textbf{Next state} & \textbf{$NR_a$} & \textbf{$NR_b$} & \textbf{$RA$} \\ \hline
$c^{(0)}_{t+1}$     & $c^{(1)}_t$     & $c^{(1)}_t$     & $c^{(0)}_t$   \\
$c^{(1)}_{t+1}$     & $c^{(1)}_t$     & $c^{(2)}_t$     & $c^{(0)}_t$   \\
$c^{(2)}_{t+1}$     & $c^{(1)}_t$     & $c_t$           & $c^{(0)}_t$   \\
\Xhline{2\arrayrulewidth}
\end{tabular}
\end{table}

\section{Optimal Battery Control  Algorithm} \label{sec:bc}

Thanks to our proposed model of  instantaneous degradation cost, {\KB we can define the MDP form for the OBC problem.} 
%Along with the additional state elements $c^{(i)}_t$ and cycle-based degradation cost, state, transition kernel and reward are updated to adapt rainflow algorithm. 
First, each state is given by
%includes three more variables $c^{(i)}_t$ which is modifed as
\begin{align}
s_t = [p_t, f_t, c_t, c^{(0)}_t, c^{(1)}_t, c^{(2)}_t],  \quad\forall t \in \ccalT \label{rain_state}
 \end{align}
which is used to determine the action $a_t$ based on the policy of interest. The transition kernel $\mathcal{P}$ now includes the updates in Table~\ref{tb_transition}, while the instantaneous reward is given by
\begin{align}
%\mathcal{R}: 
r_t (s_t, a_t) = -h^e_t - h^f_t - h^d_t, \quad \forall t \in \ccalT. \label{reward_modified}
\end{align}
%{\KBB Note that even though both instantaneous reward and total reward can be negative according to Eq.~\eqref{reward_modified}, the economic profit from the battery operation becomes positive. First, the sum of the net cost for electricity usage over the time period $\mathcal{T}$ becomes positive with the energy arbitrage. Second, the profit from participating in the FR becomes positive even with the penalty by one-time payoff of FR contract regardless of the battery control.} 

%Based on the MDP including the rainflow algorithm, the objective function of the optimization problem in 
The battery control problem now becomes to determine the best policy $\pi$ for forming the action as $a_t\sim \pi(s_t)$ with $s_t$ given in Eq.~\eqref{rain_state}. To simplify the policy search, we are particularly interested in the set of parameterized policies given by $\pi_\theta (\cdot) = \pi (\cdot; \theta)$, with parameter $\theta$ optimized through
 \begin{align}
% \begin{split}
\max_{\theta }~& \mathbb{E}_{ \pi_\theta} \left[ \sum_{t=1}^{T} \gamma^{t} r_t (s_t,~a_t) \right]. %\\
%= \max_{\pi \in \Pi}~& \mathbb{E}_{s,a \sim \pi(\cdot)} [\sum_{t=0}^{T} \gamma^{t} (-h^e_t - h^f_t - h^d_t)]. 
\label{obj_theta}
%\end{split}
\end{align}
%Note that policy $\pi$ generates the mapping between the state $s_t$ and the action $a_t$ and the objective is to find the optimla policy $\pi^{\ast}$ that maximize the total rewards. 

%If the policy $\pi$ in \eqref{obj} is parameterized by $\theta$, the problem is changed from finding optimal policy $\pi^{\ast}$ to optimal parameter $\theta^{\ast}$. As a result, the objective function becomes
%\begin{align}
%\max_{\theta} \mathbb{E}_{s,a \sim \pi(\cdot)} [\sum_{t=0}^{T} \gamma^{t} r_t]. \label{obj_theta}
%\end{align} 

To solve Eq.~\eqref{obj_theta}, we can adopt certain RL algorithms to search for the optimal parameter $\theta$; see e.g., \cite{RL}. We use the deep Q-networks (DQNs) {\cite[Ch.~9]{MDP}} here as a popular RL approach based on nonlinear neural network modeling. Accordingly, the parameter $\theta$ represents the DQN weights  to be learned, and the DQN is used to obtain the so-termed Q-network that models the MDP's \textit{action-value, or Q-function}, namely the expected total future reward under a given pair of state and action:
\begin{align}
Q(s_t, a_t) := \mathbb{E}_{\pi_\theta} \left[\sum_{\tau=t}^{T} \gamma^{(\tau-t)} r_\tau (s_\tau, a_\tau) \Big \vert s_t, a_t \right].
\end{align}
For the optimal Q-function, the Bellman optimality condition \cite{RL}  states that: 
\begin{align}
%\begin{split}
&Q^{\ast}(s_t, a_t) = r_t (s_t,a_t) + \gamma \mathbb{E}_{s_{t+1}} \left[\max_{a_{t+1}} Q^{\ast}(s_{t+1}, a_{t+1} ) \Big \vert s_t, a_t\right]. \label{eq_bellman}
%\end{split}
\end{align}

To find the optimal $Q^\ast$, we parameterize {\KB the action-value Q-function} using $\theta$ as the NN weights, as denoted by $Q(s_t, a_t; \theta)$. The Bellman optimality in Eq.~\eqref{eq_bellman} can be used to develop iterative {\textit{gradient descent}} updates to obtain the best $\theta$. {At each update,} the Q-network on the right-hand side of Eq.~\eqref{eq_bellman} is kept constant as the \textit{target network}, whereas the other one is varied to minimize the difference between both sides. Letting {$\theta'$} denote the latest NN weights, we design the loss function for DQN training as the expected squared difference:
%{\hao (the following notation has issue. $\theta_{n-1}$ (why n-1?) should be inside Q input??, plz be specific with over which the expectation, Q-fn already has expectation as shown in (16))}
\begin{align}
%\begin{split}
\ccalL(\theta)=\mathbb{E}_{\{s_t, a_t,  s_{t+1}\}}&\Big[\big(r_t + \gamma\max_{a_{t+1}}Q(s_{t+1}, a_{t+1};\theta')\nonumber\\
&-Q(s_t, a_t;\theta) \big)^2 \Big].  \label{loss_function}
%\end{split}
\end{align}
%where $\theta_{n-1}$ is used to fix the target Q-function in the right-hand side of \eqref{eq_bellman}.
%{\KB Note that $\theta_n$ are the parameter of the neural network at the $n$th iteration, and the parameter in the previous iteration $\theta_{n-1}$ is held fixed when optimizing the loss function.} 
%This loss function captures the expected squared differences between both sides of \eqref{eq_bellman}, with one of the Q-functions fixed by the last iterate $\theta_{n-1}$. Basically, the right-hand side of \eqref{eq_bellman} is a fixed target for the update of iteration $n$. 
%{\hao (need to mention how $\theta_{n-1}$ is used right away, then elaborate more)}. %Specifically, the left side is computed as the approximated Q-value with weight parameter $\theta_{n}$ whereas the right side is computed with the reward $r_t$ and maximum of the approximated Q-value with weight parameter $\theta_{n-1}$. Note that the right side with paratmer $\theta_{n-1}$ works as the target for iteration $n$ to update $\theta_{n}$ at the next iteration.
To minimize $\ccalL(\theta)$, one can need to compute its gradient over the parameter $\theta$ given by 
\begin{align}
%\begin{split}
\nabla_{\theta} \ccalL(\theta)
=&\mathbb{E}_{\{s_t, a_t,  s_{t+1}\}}\Big[-2\Big(r_t + \gamma \max_{a_{t+1}}Q(s_{t+1}, a_{t+1};\theta') \nonumber\\
&-Q(s_t, a_t;\theta)\Big) \nabla_{\theta} Q(s_t, a_t; \theta)\Big]. \label{loss_gradient}
%\end{split}
\end{align}
 %{\hao (plz add details and first introduce algorithm)}
Each gradient-based update relies on the estimate from sampling the trajectory such that the expectation in Eq.~\eqref{loss_gradient} is replaced by the sample average. {To this end, the action $a_{t}$ is sampled for given state $s_{t}$ based on $\theta'$ as  $a^{\ast}_t = \argmax_{a_t} Q(s_t, a_t; \theta')$, $\forall t$.} To ensure adequate exploration of the state space, the $\epsilon$-greedy method  \cite[Ch.~2]{MDP} can be used to randomize the action by selecting $a_t^\ast$ with probability $(1-\epsilon)$ at every time. The value of $\epsilon$ would decrease as the DQN updates continue, typically at an exponential decreasing rate $\kappa\in (0,1)$. This method can improve the exploration process at the beginning phase while eventually picking the optimal actions to attain convergence. 
%{\hao (plz show the equation of estimating (19) from trajectory)}
	
{\KB To improve the efficiency and stability of DQN implementation, we introduce two additional techniques. First, we implement the \textit{experience replay} method \cite{experience_replay} to efficiently use the past samples by storing all the past samples in the memory $\mathfrak{D}:=\{(s_t, a_t, s_{t+1}, r_t)\}$ along the trajectory.} {When computing the loss function Eq.~\eqref{loss_function}, a subset of samples denoted by mini-batch $\mathfrak{J}$ is randomly picked from $\mathfrak{D}$ and used as the samples for gradient estimation. This method can improve the training efficiency by selectively reusing past samples.}
{In addition, we advocate the \textit{fixed target network} approach \cite{fixed_target_network} by keeping the target network parameter fixed for several updates. To this end, let $\theta^{-}$ denote the target network parameter, which is only updated once every $N_o$ iterations.} This technique could mitigate any potential instability issue by changing the DQN target weights less frequently. By adopting \textit{experience replay} method and \textit{fixed target network} approaches, we obtain the estimates of loss function and its gradient as %{\hao (this needs to be fixed after fixing (19) )}
\begin{align}
%\hat{\ccalL}(\theta_n)&=(1/|\mathfrak{J}|) \sum_{t \in \sim\mathfrak{J}} \big[(r_t + \gamma\max_{a_{t+1}}Q(s_{t+1}, a_{t+1};\theta^{-}))\nonumber\\
%&\hspace*{20mm} -Q(s_t, a_t;\theta_n)\big]^2,  \label{loss_function_final}\\
\nabla_{\theta} \hat{\ccalL}(\theta)
&=(1/|\mathfrak{J}|) \sum_{t \in \mathfrak{J}}\Big[-2 \Big(r_t + \gamma \max_{a_{t+1}}Q(s_{t+1}, a_{t+1};\theta^{-}) \nonumber\\
&\hspace*{5mm}-Q(s_t, a_t;\theta)\Big) \nabla_{\theta} Q(s_t, a_t; \theta)\Big].
\label{loss_gradient_final}
\end{align}

The detailed algorithmic steps for DQN-based OBC algorithm are tabulated in Algorithm \ref{alg:dqn}. %{\hao (need some discussions on convergence, note that usually RL doesn't go to the max iteration $N$)} 
{As mentioned earlier, the state variables $p_t$ and $f_t$ are not action dependent. Thus, their transitions are obtained from the profiles given as the algorithm input, such as real data provided by the market operators.}
{\KB For the convergence of DQN algorithms, the total number of episodes $N$ is typically chosen to be large enough in practice. For each episode $n$, there are total $T$ samples from $t=1$ to $t=T$. Note that Algorithm \ref{alg:dqn} can be used to search for the best policy under the linearized degradation cost as well, by using this simpler degradation cost model in Eq.~\eqref{linear_deg}. The ensuing section will compare these two degradation models numerically.}
%{\hao (why need (20)? can we just explain it in words?)}
 %The deep Q-learning algorithm is indicated in \textbf{Algorithm 1}.
% Have seen the above part
\begin{algorithm}[t]
\SetAlgoLined
\caption{DQN-based Optimal Battery Control}
\label{alg:dqn}
\DontPrintSemicolon
{\bf Hyperparameters:}  discount factor $\gamma = 1$, learning rate $\eta>0$, $\epsilon$-greedy coefficient $\kappa \in (0,1)$, mini-batch size $|\mathfrak{J}|$, target network update interval $N_o$, and maximum number of {\KB episodes} $N$. \;
{\bf Input:} training profiles of prices and FR signals with the exploration time horizon $T$. \;
{\bf Initialize:} the $\epsilon$-greedy probability $\epsilon \in (0,1)$, replay memory $\mathfrak{D}=\emptyset$, {initial action-value function $Q(s, a ; \theta')$ with a random $\theta'$ and  the target network parameter $\theta^- =\theta'$ at episode $n=0$.} \;
\While{$n \leq N$}{
%{	\hao perhaps use subscript $(n)$ as iteration index} \;
%observe $s_{n,0}$ at the initial $t=0$ in $n$ th iteration, sample \;
%\For{profile $1$ to profile $N_p$}{
\For{t=1, $\cdots$ ,T}{
Select a random action $a_{t}$ with probability $\epsilon$; otherwise,   use the action $a_{t}^* = \argmax_{a_t} Q(s_{t}, a_t; \theta')$.\;
Implement the action $a_{t}$ to obtain the ensuing state $s_{t+1}$ based on {the transitions of both Eq.~\eqref{transition_cc} and Table \ref{tb_transition}}, and by using the input profiles of $p_{t+}$ and $f_{t+1}$. \;
Compute the instantaneous reward $r_{t}$ in Eq.~\eqref{reward_modified}.\;
Store the tuple $(s_{t}, a_{t}, s_{t+1}, r_{t})$ in $\mathfrak{D}$.\;
Select a random mini-batch $\mathfrak{J}$ with size $|\mathfrak{J}|$ from $\mathfrak{D}$.\;
Compute the gradient estimate using Eq.~\eqref{loss_gradient_final}.\;
{Update the parameter $\theta' \leftarrow \theta' - \eta  \nabla \hat \ccalL(\theta')$.\;
\If{$t/N_o$ is an integer}{Update the target network parameter $\theta^{-} \leftarrow \theta'$.}}
Update $\epsilon = \kappa \epsilon$.
}
Update the episode number $n\leftarrow n+1$. }
\end{algorithm}

%%%%%%%%%%%%%%%%%%%%%%%%%%%%%%%%%%%%%%%%%%%%%%%%%%%%%%%%%%%%%%%%%%%%%
 %%
 %%      Section: Case Study %%
 %%%%%%%%%%%%%%%%%%%%%%%%%%%%%%%%%%%%%%%%%%%%%%%%%%%%%%%%%%%%%%%%%%%%
%%

\section{Numerical Tests} \label{sec:CS}

%\subsection{System Description}
{\KB We have compared the proposed RL-based battery control algorithm under cycle-based degradation cost with the linearized approximation one \cite{4-6}.} Actual data of electricity market prices and FR signals have been used, respectively from the ERCOT's market data depository \cite{ERCOT} and PJM's ancillary service datasets \cite{PJM}. Each time instance corresponds to a 5-minute interval. The FR signal is normalized to indicate either maximum charging or discharging for the battery. {\KBA As mentioned in Remark \ref{rmk:frsignal}, the fast FR signal at 2-second rates is averaged over a 10-second interval for the training phase, while the original data rate is maintained for the testing phase. We have used a 200kWh-capacity battery with (dis)charging rate of 120kW and minimum SoC of 20kWh, which takes 90 minutes to fully (dis)charge. The multiple discrete action space is adopted with overall 11 actions, as $\mathcal{A} = \{-1, -0.8, \cdots, 0.8, 1 \}$. The parameters associated with battery degradation are set to $\alpha_d = 4.5 \times 10^{-3}$ and $\beta = 1.3$, as used in \cite{4-3}.}

%slots in the sampling and thus reduce the time to converge in the training. %
%Examples profiles of normalized electricity market price and FR signal are depicted in Fig.~\ref{fig_inputdata}. 
%\begin{figure}[t]
%\centering
%\includegraphics[width=.75\linewidth,height=70mm]{fig_inputdata2.png}
%\caption{Example profiles of normalized state variables: (a) electricity market price $p_t$ and (b) FR signal $f_t$. 
	%{\hao (may not need everything)}) 
%}
%\label{fig_inputdata}
%\end{figure}

The DQN \textbf{Algorithm 1} has been implemented in Python with the  popular NN toolboxes Tensorflow and Keras \cite{Keras}. Table~\ref{tb_training_setting} lists the parameter settings for  the DQN training, which uses {7 daily profiles of $\{p_t,f_t\}$.} {\KB There are a total of $T=8,640$ time instances for each exploration episode. Upon the convergence of Q-network, it is used for determining the optimal (dis)charging actions for each 2-second interval of 60 days of testing data, while each testing trajectory having 43,200 time instances.}

\begin{table}[t]
    \centering
    \caption{\KBC Parameter settings for DQN training}
    \label{tb_training_setting}
    \begin{tabular}{ c  c }
        \toprule
        \textbf{Parameter}					&	\textbf{Value}\\ \hline
        Number of hidden layers     & 2\\
        Number of nodes        &  [128, 32]\\
        Activation function        & ReLU\\
        Learning rate			&0.001\\
        Optimizer				& Adam\\
        Epsilon ($\epsilon$)		&0.001\\
        Batch size ($J$)				& 256\\
        Maximum number of episodes ($N$) 		& 2000\\
        Number of daily profiles	& 7\\
        \bottomrule
    \end{tabular}
\end{table}

\begin{figure}[t]
\centering
\begin{subfigure}{\linewidth}
\includegraphics[width=\linewidth]{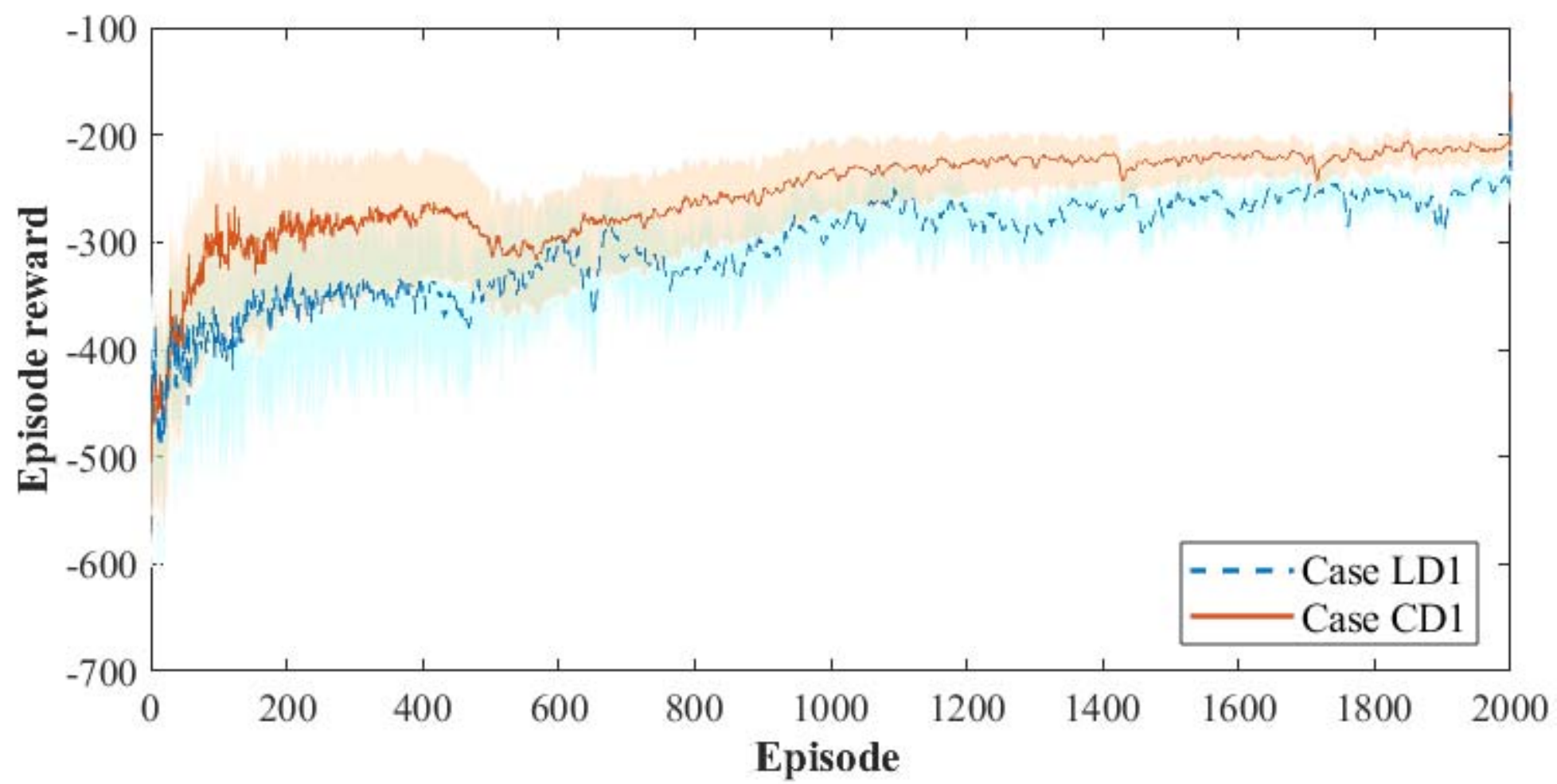}
\caption{}
\end{subfigure}
\begin{subfigure}{\linewidth}
\includegraphics[width=\linewidth]{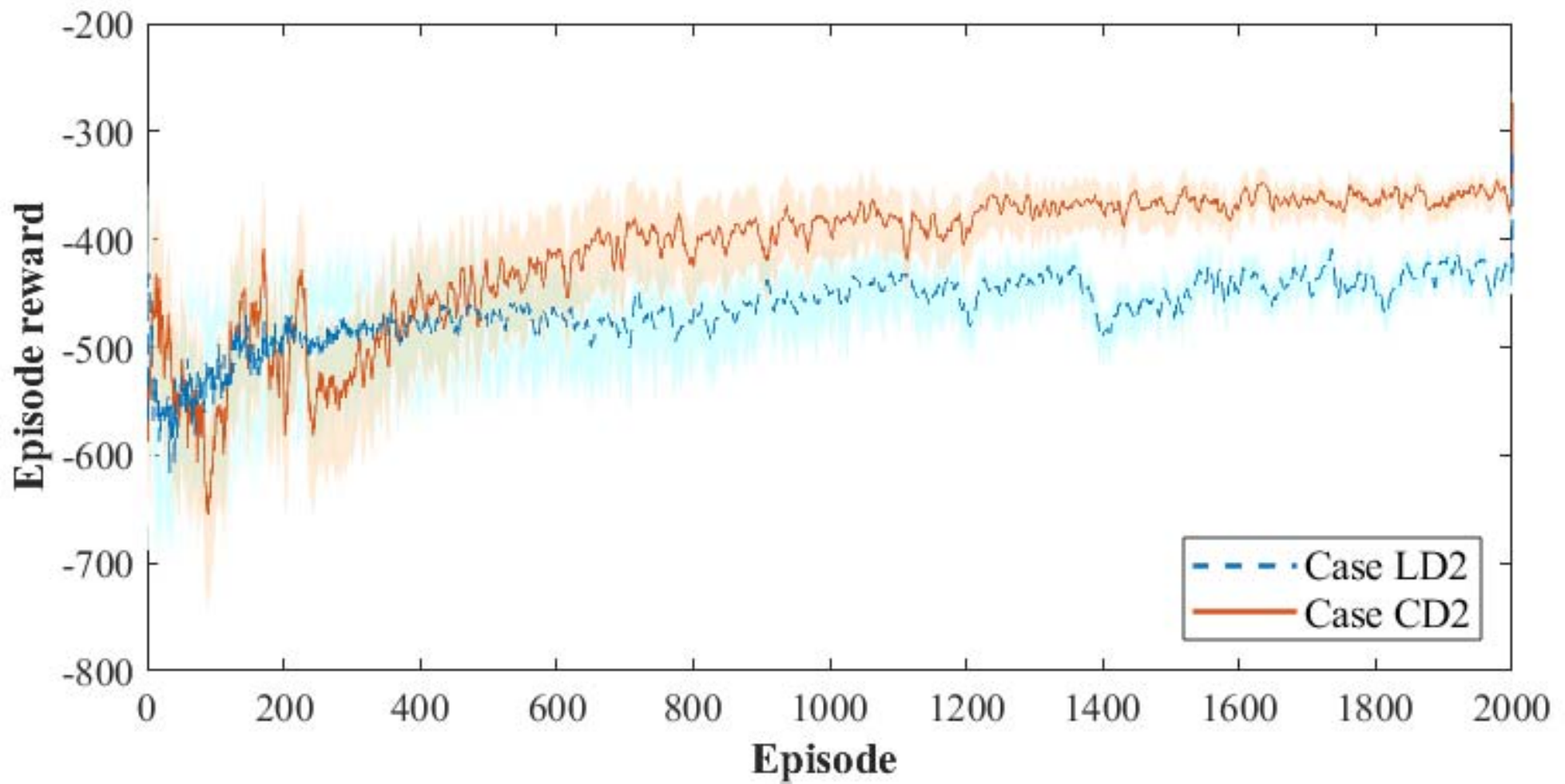}
\caption{}
\end{subfigure}
\caption{Comparisons of the {\KB total} reward trajectory between (a) cases LD1 and CD1 ($\alpha_d$) and (b) case LD2 and CD2 ($2\alpha_d$)} %{\hao (don't see the value of including these figures)}}
	\label{fig_convergence}
\end{figure}

\textbf{Training Comparisons.} To compare the proposed cycle-based degradation cost with the linear one (denoted by CD and LD, respectively) in terms of  battery control performance, We have considered two levels of degradation coefficients, at $\alpha_d$ and $2\alpha_d$, respectively. Fig.~\ref{fig_convergence} illustrates the training comparisons of the actual episode rewards (all based on cycle-based degradation) for all the three cases. Clearly, all the DQN iterations are convergent as the {\KB total} reward trajectories tend to be non-decreasing till reaching the highest values. Moreover, while the CD cases using our proposed \textit{instantaneous cycle-based degradation} modeling outperform the LD counterparts, especially at larger degradation cost. This comparison validates the advantages of the proposed degradation model in terms of accurately representing the battery cost and thus leading to effective control policies. 

\begin{figure}[t]
\centering
\begin{subfigure}{\linewidth}
\includegraphics[width=\linewidth]{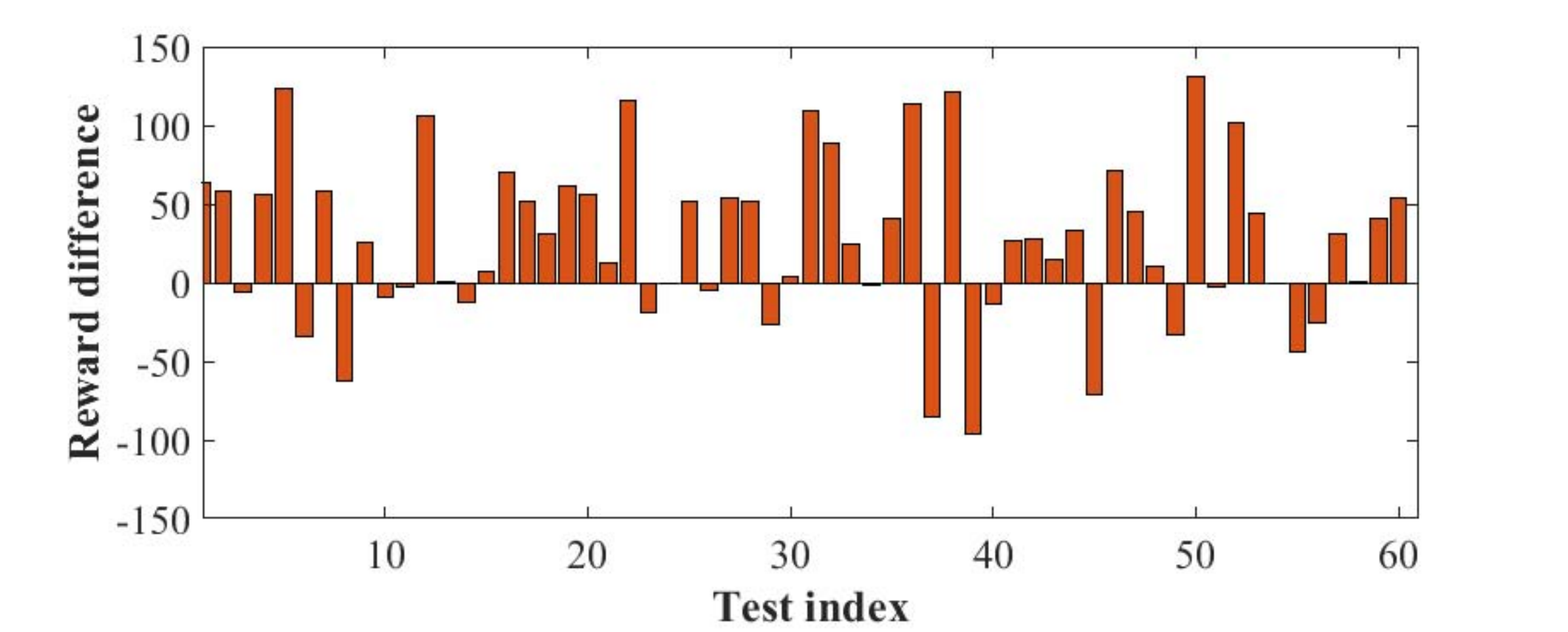}
\caption{}
\end{subfigure}
\begin{subfigure}{\linewidth}
\includegraphics[width=\linewidth]{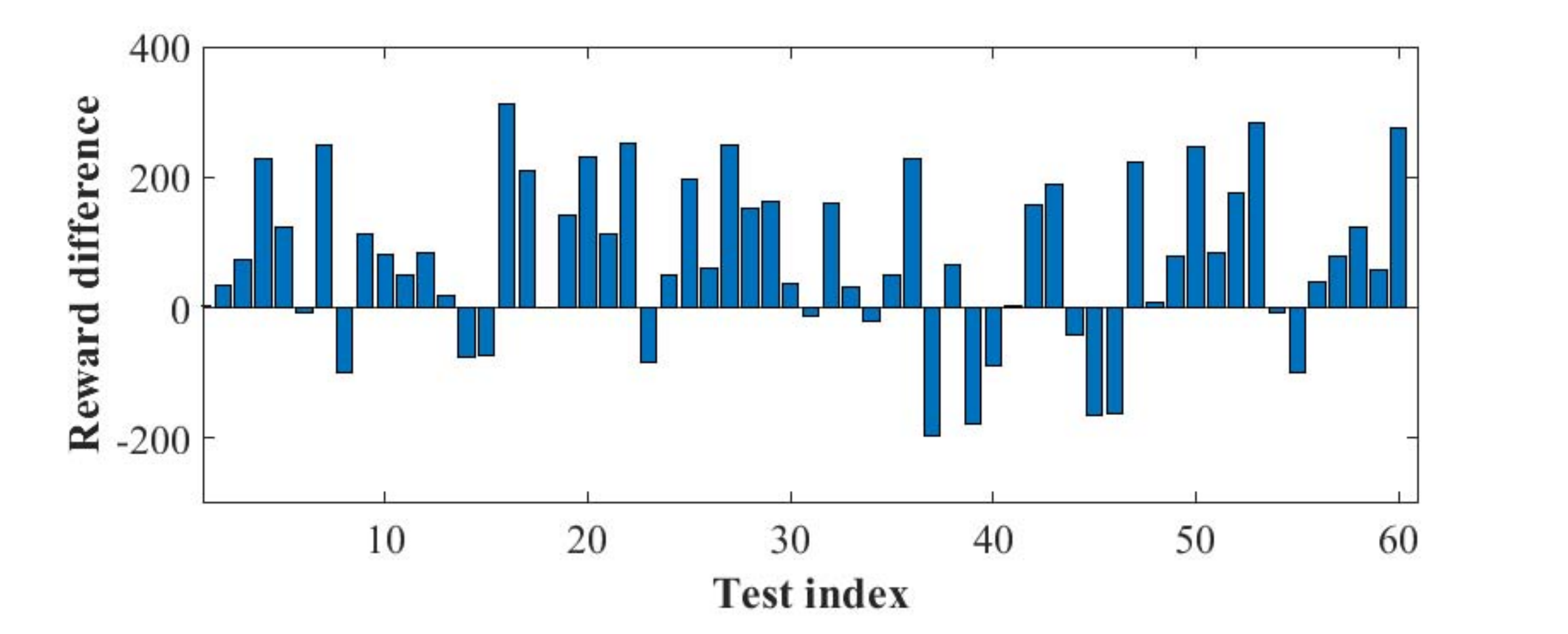}
\caption{}
\end{subfigure}
\caption{The total reward differences (positive difference indicating higher reward for CD) between (a) cases LD1 and CD1 ($\alpha_d$) as well as (b) case LD2 and CD2 ($2\alpha_d$)}
\label{fig_test}
\end{figure}

\begin{figure}[t]
\centering
\begin{subfigure}{\linewidth}
\includegraphics[width=\linewidth]{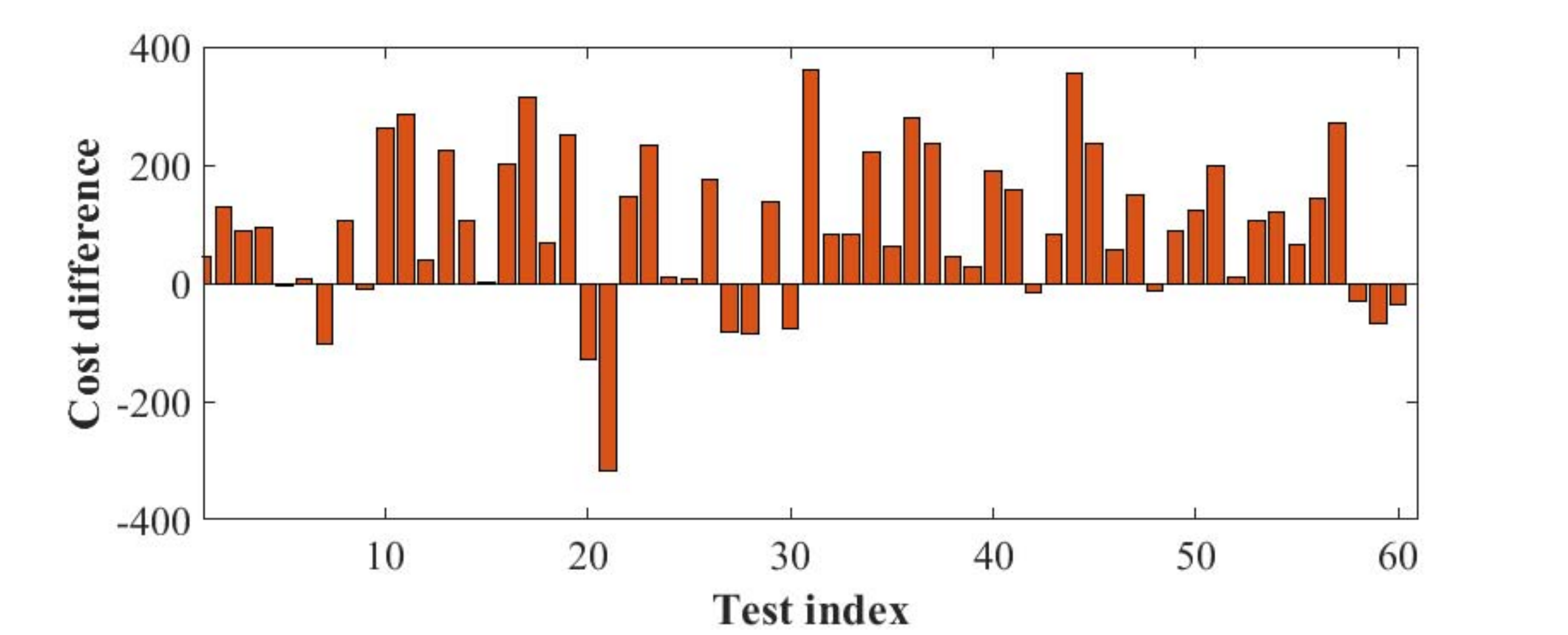}
\caption{}
\end{subfigure}
\begin{subfigure}{\linewidth}
\includegraphics[width=\linewidth]{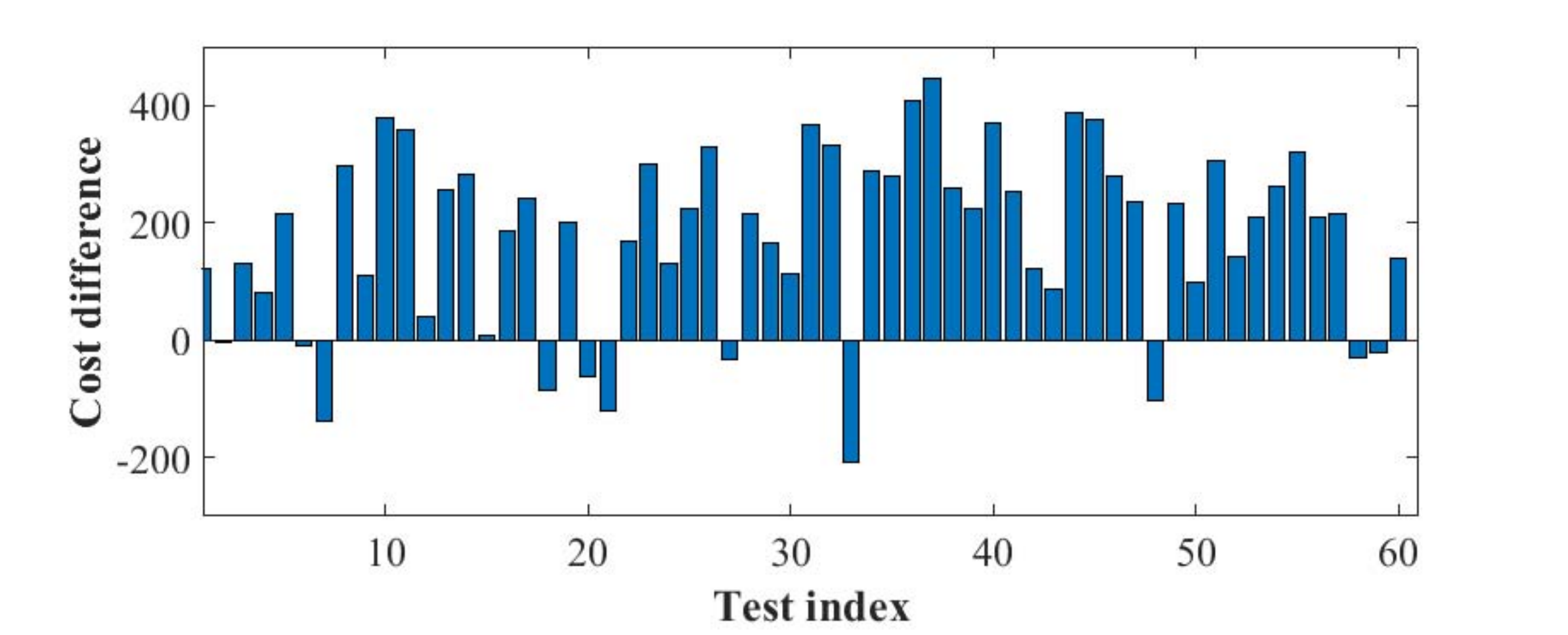}
\caption{}
\end{subfigure}
\caption{The battery degradation cost differences (positive differences indicating lower cost for CD) between (a) cases LD1 and CD1 ($\alpha_d$) as well as (b) case LD2 and CD2 ($2\alpha_d$)}
\label{fig_deg}
\end{figure}

\textbf{Testing Comparisons.}
We have further compared the testing performance of the learned Q-networks using both CD and LD based models. Fig.~\ref{fig_test} plots the total reward differences between the CD and LD solutions (positive differences indicating higher reward for CD) for each test trajectory under the two levels of degradation coefficients. The proposed CD based control leads to higher total reward for at least 73.33\% or 81.67\% of test scenarios, respectively for the two $\alpha_d$ levels. This result confirms the earlier observations in training phase that proposed solution is more attractive for larger degradation cost. {\KBB Table~\ref{tb_diff_average} indicates the total mean, maximum, minimum values and average values of the cases when CD has better reward than LD, and vice versa. As shown in the table, the overall mean value increases as the degradation coefficient increases and the battery degradation cost affects more in total cost accordingly. In addition, the maximum, minimum and mean values show that even though there are some cases that LD show better performance than CD, the number of these cases is very small.}
Similar comparison is also observed in differences of battery degradation performance only (again, positive differences indicating lower degradation cost for CD), as illustrated by Fig.~\ref{fig_deg}. Clearly, the proposed CD solutions overwhelmingly improve the battery degradation performance and accordingly the total reward, as compared to the existing LD-based approximation. {\KBA In addition, Fig.~\ref{fig_test} and  Fig.~\ref{fig_deg} share very similar pattern, which implies that the increase in total reward is mostly caused by the decrease in the battery degradation cost and has least impacts on the decrease in the rewards regarding the net cost for electricity usage or frequency regulation penalty.} %Similar comparison  is also observed in differences of battery degradation performance only (again, positive differences indicating lower degradation cost for CD), as illustrated by Fig.~\ref{fig_deg}. Clearly, the proposed CD solutions overwhelmingly improve the battery degradation performance and accordingly the total reward, as compared to the existing LD-based approximation. {\KB In addtion, Fig.~\ref{fig_test} and  Fig.~\ref{fig_deg} share the similar patterns in most tests, which implies that the increase in total reward is mostly caused by the decrease in the battery degaradation cost and has least impacts on the decrease in the rewards regarding the net cost for electricity usage or frequency regulation penalty.}

\begin{table}[t]
	\centering
	\caption{\KBC Reward differences between CD and LD}
	\label{tb_diff_average}
	\begin{tabular}{ c | c  c  c|  c  c }
		\Xhline{2\arrayrulewidth}
		\makecell{ \textbf{Cases} \\ \textbf{(A-B)}} &\textbf{Mean} &\textbf{Max} &\textbf{Min}  &\makecell{ \textbf{Mean} \\ \textbf{(A > B)}} &\makecell{ \textbf{Mean} \\ \textbf{(A < B)}}\\ \hline
		CD1-LD1    	&84.45 & 133.09 & -97.04 & 111.38 & -72.12\\
		CD2-LD2        	&176.51 & 315.15 & -198.85 &172.79 & -134.19\\
		\Xhline{2\arrayrulewidth}
	\end{tabular}
\end{table}

%Fig.~\ref{fig_test} indicates that in most of the cases, case CD shows bigger reward compared to case LD. The percentage that case CD has bigger reward than case LD among 60 tests are  It implies that we can find more optimal policy when modeling the battery degradation cost using \textit{instantaneous cycle-based degradation cost}. 
%\begin{figure}[t]
%\centering
%\includegraphics[width=\linewidth]{fig_deg.png}
%\caption{The battery degradation differences (positive differences indicating lower degradation cost for CD) between (a) cases LD1 and CD1 ($0.5\alpha_d$), (b) case LD2 and CD2 ($\alpha_d$), as well as (c) case LD3 and CD3 ($2\alpha_d$).}
%\label{fig_deg}
%\end{figure}

\begin{figure}[t]
\centering
\begin{subfigure}{\linewidth}
\includegraphics[width=\linewidth]{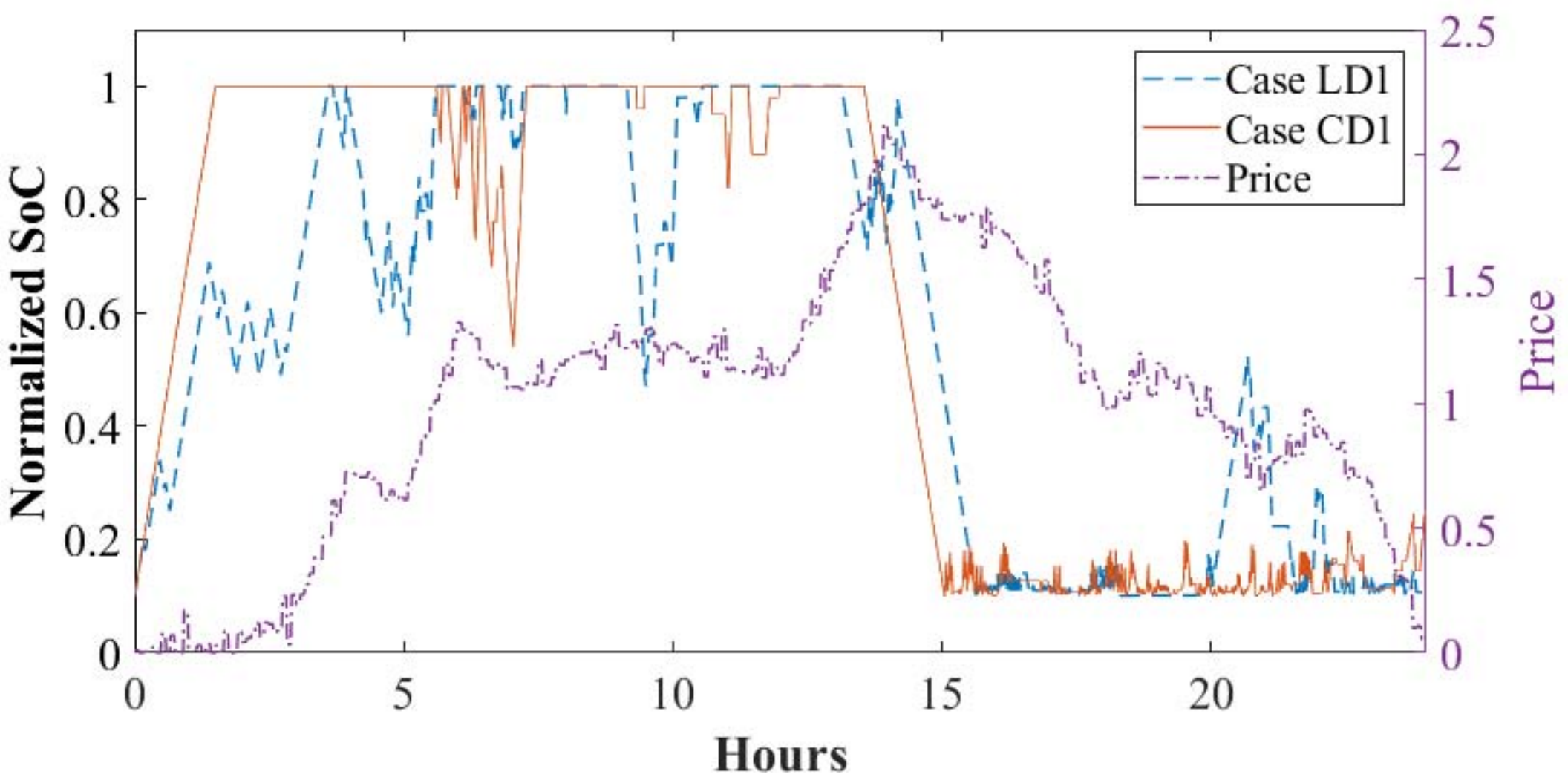}
\caption{}
\end{subfigure}
\begin{subfigure}{\linewidth}
\includegraphics[width=\linewidth]{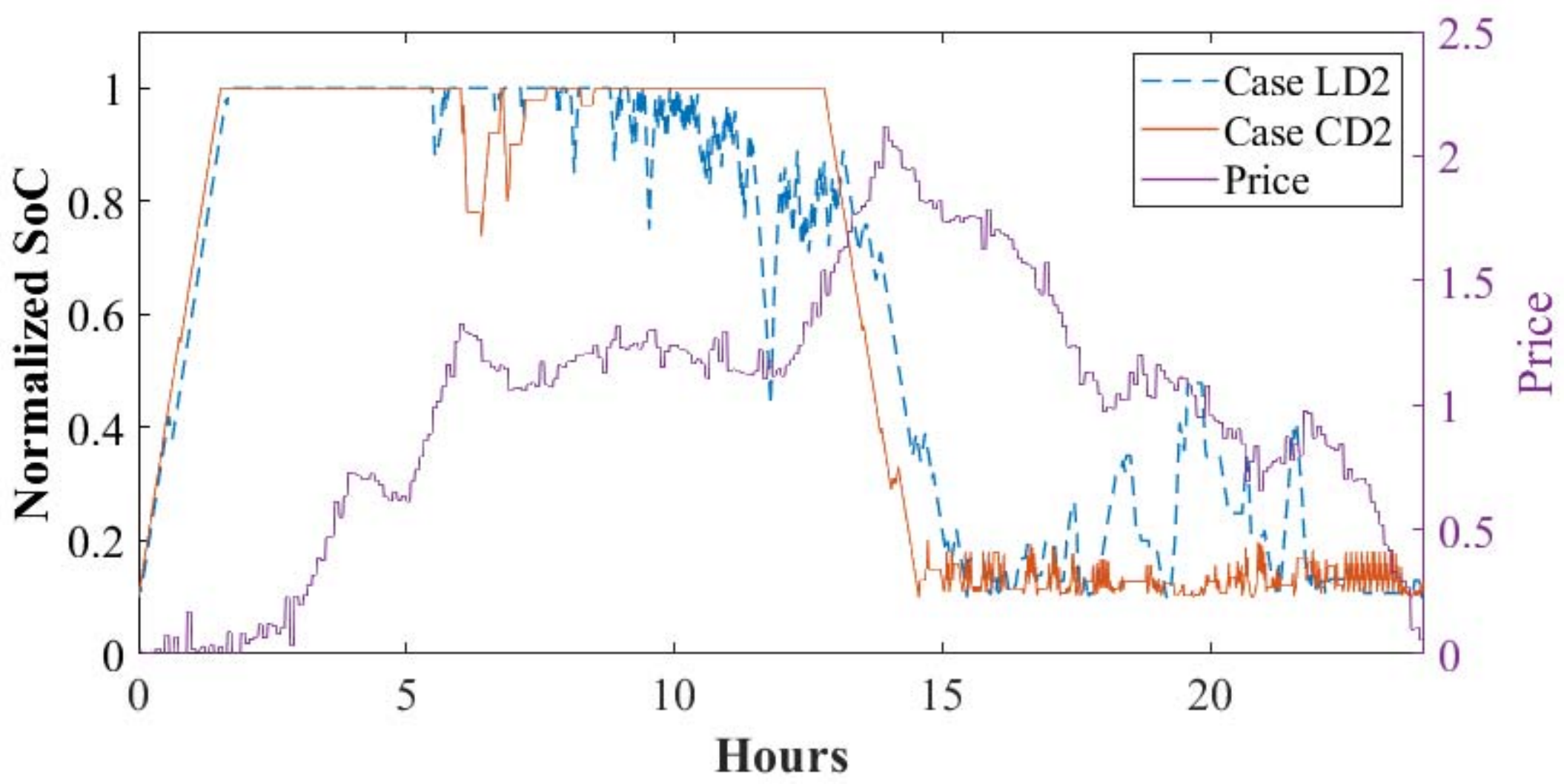}
\caption{}
\end{subfigure}
\caption{Comparison of selected SoC trajectories in testing between (a) cases LD1 and CD1 ($\alpha_d$) and (b) cases LD2 and CD2 ($2\alpha_d$)}
\label{fig_soc}
\end{figure}

Fig.~\ref{fig_soc} plots the selected testing SoC trajectory {\KB along with the electricity price} to better illustrate the improvement of the proposed CD-based policy, corresponding to the two choices of degradation parameter {\KB ($\alpha_d$ and $2\alpha_d$)}. {\KBB Clearly, both trajectories show that the CD-based policy leads to less number of cycles with long depth as compared to the LD one, which reduces the overall degradation cost especially for hours between $[3,~13]$.} {\KBB In addition, during high-price hours in $[12,~15]$, the CD trajectory has one smooth and long discharging cycle and this pattern is amendable to mitigating battery degradation.} In contrast, the LD one has frequent, noticeable fluctuations during this period. Because of the linearized approximation, the LD-based policy leads to eight more noticeable (dis)charging cycles of considerate depth than the CD one. {This speaks for the capability of the proposed CD model in effectively removing some unnecessary cycles of moderate depth, thanks to the accurate representation of rainflow-based degradation.} In addition, in the post-peak hours {\KB $[15,~20]$}, the LD based policy produces a couple of cycles of moderate depth which are not very profitable. The proposed CD based policy is able to successfully remove these nonprofitable cycles and does not lead to any considerate cycles. 

\begin{table}[t]
	\centering
	\caption{{\KBD Degradation comparisons between CD and LD}}
	\label{tb:deg_comp}
	\begin{tabular}{ c | c c | c c }
		\Xhline{2\arrayrulewidth}
		\textbf{Degradation factor} & \textbf{LD1} & \textbf{CD1} & \textbf{LD2} & \textbf{CD2}\\ \hline
		High C-rates & 0.0563 & 0.0342 & 0.0463 & 0.0369 \\
		SoC stress & 0.0141 & 0.0108 & 0.0182 & 0.0107\\ \Xhline{2\arrayrulewidth}
	\end{tabular}
\end{table}	

Interestingly, by mitigating cycle-based degradation the proposed CD approach can potentially contribute to the improvement of other degradation factors too. {\KBD Table~\ref{tb:deg_comp} compares the proposed CD with the LD approach on the degradation related to high C-rates \cite{wang} and SoC stress \cite{4-4}, both of which have been numerically improved by the CD-based policies.} The high C-rate based degradation depends on the total DoD summed over all cycles of the trajectory. Intuitively, a concise list of smooth and long (dis)charging cycles attained by CD-based policy can reduce both the number of cycles and their DoD, thus beneficial for the high C-rate metric. Similarly, as CD-based policy has also been observed to remove unnecessary cycles in the post-peak hours, {the average SoC level decreases which relieves the SoC stress.} These intuitions corroborate the claim in Section \ref{sec:cost} that the cycle-based DoD stress model is most relevant for the fast battery control problem.

To sum up, the numerical results have validated the performance improvement attained by the proposed \textit{instantaneous} cycle-based degradation model, by exactly representing the rainflow conditions. The proposed approach effectively leads to battery control trajectories that reduce unnecessary fluctuations or improve the overall economical profits.

%%%%%%%%%%%%%%%%%%%%%%%%%%%%%%%%%%%%%%%%%%%%%%%%%%%%%%%%%%%%%%%%%%%%%
 %%
 %%      Section: conclusions %%
 %%%%%%%%%%%%%%%%%%%%%%%%%%%%%%%%%%%%%%%%%%%%%%%%%%%%%%%%%%%%%%%%%%%%
%%

\section{Conclusion and Future Work} \label{sec:CON}

The paper proposes an accurate model of cycle-based degradation cost in order to allow for efficient battery control designs using reinforcement learning (RL). In order to model the degradation which depends on the full cycle, we introduce additional state variables to judiciously keep track of important switching points of SoC trajectory for effectively identifying (dis)charging cycles. This way, the actual degradation cost is separated into instantaneous terms along with other operation costs such as the net cost for electricity usage and FR penalty, such that  powerful DQN based RL algorithms are readily applicable. Numerical tests confirm the effectiveness of proposed cycle-based degradation model and demonstrate the  performance improvements in effectively mitigating battery degradation over existing linearized approximation approach. 

{Exciting future research directions open up on expanding the battery operations to support other ancillary services such as peak-shaving. The key question will be how to model the peak-shaving cost as instantaneous reward. In addition, the power network constraints such as voltage limit are of high interest in practice and will be considered as well.}
%{\hao (need a paragraph on future research)}

%%%%%%%%%%%%%%%%%%%%%%%%%%%%%%%%%%%%%%%%%%%%%%%%%%%%%%%%%%%%%%%%%%%%%%
%                                                                    %
%       Appendix %
%
%%%%%%%%%%%%%%%%%%%%%%%%%%%%%%%%%%%%%%%%%%%%%%%%%%%%%%%
%\numberwithin{equation}{subsection}
{\KBC 
\appendix
%\subsection{Proof of Proposition 1}

The goal is to show that Eq.~\eqref{inst_deg_cost} represents the exact incremental degradation cost from time $t$ to $(t+1)$. With the current SoC denoted as $c_t$ at time $t$,  the current cycle depth equals to $d_t = c_t-c_t^{(2)}$  starting from the SP  $c_t^{(2)}$ in time $t_0$. Thus, the resultant degradation cost is given by
\begin{align*}
	\Phi(d_t) = \alpha_d e^{\beta d_t} = \alpha_d e^{\beta \vert c_t - c_t^{(2)} \vert}
\end{align*}
with the absolute difference capturing cycle depth. When transitioning to time $(t+1)$, the incremental difference in degradation cost due to action $b_t$ becomes
\begin{align*}
	\Delta_t &= \Phi(d_{t+1}) - \Phi(d_t) = \Phi(d_t + b_t) - \Phi(d_t) \nonumber \\
	&= \alpha_d e^{\beta \lvert c_t+b_t-c^{(2)}_t\rvert}-\alpha_d e^{\beta \lvert c_t-c^{(2)}_t \rvert}.
\end{align*}
Therefore, summing up these differences leads to
\begin{align}
	\Phi(d_{t+1}) &= \Phi(d_t) + \Delta_t = \textstyle \sum_{\tau=t_0+1}^{t+1} \Delta_{\tau} \label{phi_sum}
\end{align}
with the degradation cost at time $t_0$ initialized by zero.
Notably, there are two types of scenarios when considering this summation at time $(t+1)$, as detailed here.
 
% \begin{enumerate}
\noindent {\textit{(1) Case $NR_a$ and case $NR_b$}}: In either case, rainflow condition is not satisfied as depicted in Fig.~\ref{fig_no_rainflow_cases}, and thus Eq.~\eqref{phi_sum} accumulates the total degradation thus far for this cycle.

\noindent \textit{(2) Case $RA$}: 
	Without loss of generality, consider the $RA$ cycle $A-K-B-C-L-D$ as shown in Fig.~\ref{fig_rainflow_example}. As the rainflow condition is satisfied at point $L$, the total degradation cost of this cycle  equals to 
%	\begin{align*}
		$\bar{\Phi} = \Phi(d_0) + 2\Phi(d_1)$.
%	\end{align*}
	Following from Eq.~\eqref{phi_sum}, the degradation cost from point $A$ to $L$ is obtained by 
	\begin{align*}
		{\Phi}^{A\rightarrow L} = \Phi(d_K+d_1) +2\Phi(d_1).
	\end{align*}
 Similarly, from point $L$ to $D$, the cycle continues on as $A-B-L-D$ due to the satisfaction of rainflow condition, and Eq.~\eqref{phi_sum} leads to an additional degradation cost as
	\begin{align*}
		{\Phi}^{L\rightarrow D}= \Phi(d_0) - \Phi(d_K+d_1).
	\end{align*}
	Together, the total degradation cost for the  cycle $A-K-B-C-L-D$ equals to $\bar{\Phi}$, which completes the proof for Proposition \ref{prop:action}. } \hfill \qedsymbol 
%\end{enumerate}

%%%%%%%
%%%%%%%%%%%%%%%%%%%%%%%%%%%%%%%%%%%%%%%%%%%%%%%%%%%%%%%%%%%%%%%
%                                                                    %
%       Bibliography %
%
%%%%%%%%%%%%%%%%%%%%%%%%%%%%%%%%%%%%%%%%%%%%%%%%%%%%%%%
%\section*{References}
%\vspace{-10pt}
\bibliography{bibliography}
\bibliographystyle{IEEEtran}
\itemsep2pt

%\newpage
\section*{Biographies}
\vspace{-2cm}
\begin{IEEEbiography} [{\includegraphics[width=1in,height=1.25in,keepaspectratio]{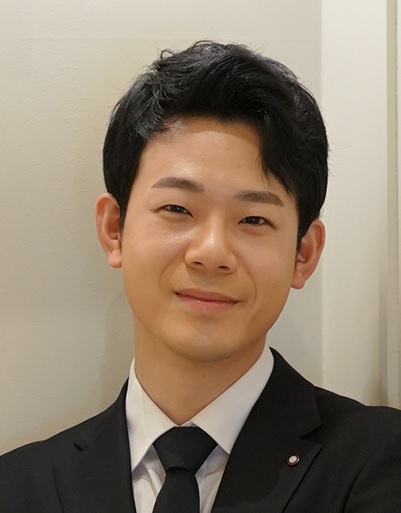}}]
	{Kyung-bin Kwon} (S'21) received the B.S. and M.S. degrees in electrical and computer engineering from Seoul National University, Seoul, South Korea, in 2012 and 2014, respectively. He is currently working towards the Ph.D. degree at The University of Texas at Austin, Austin, USA. His current research focuses on reinforcement learning applications in power system and data-driven control of networked system.
\end{IEEEbiography}
\vspace{-2cm}
\begin{IEEEbiography}[{\includegraphics[width=1in,
		height=1.25in,keepaspectratio]{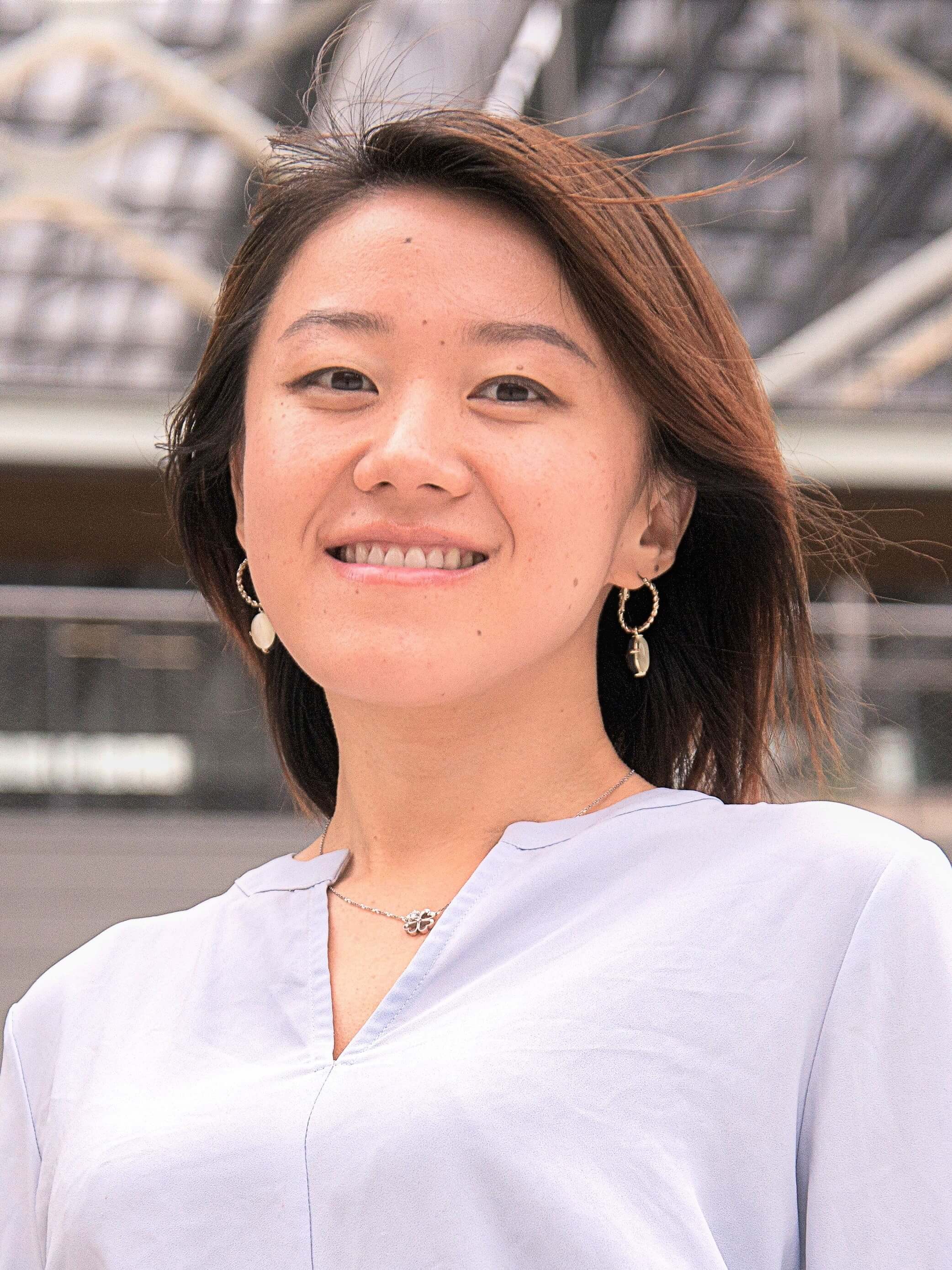}}]{Hao Zhu} (M'12--SM'19) is an Associate Professor of Electrical and Computer Engineering (ECE) at The University of Texas at Austin. She received the B.S. degree from Tsinghua University in 2006, and the M.Sc. and Ph.D. degrees from the University of Minnesota in 2009 and 2012. From 2012 to 2017, she was a Postdoctoral Research Associate and then an Assistant Professor of ECE at the University of Illinois at Urbana-Champaign. Her research focus is on developing innovative algorithmic solutions for problems related to learning and optimization for future energy systems. Her current interest includes physics-aware and risk-aware machine learning for power system operations, and energy management system design under the cyber-physical coupling. She is a recipient of the NSF CAREER Award and an invited attendee to the US NAE Frontier of Engr.~(USFOE) Symposium, and also the faculty advisor for three Best Student Papers awarded at the North American Power Symposium. She is currently an Editor of \textit{IEEE Trans.~on Smart Grid} and \textit{IEEE Trans.~on Signal Processing}. 
\end{IEEEbiography}

%\section*{Acknowledgments}
%This work has been supported by NSF Grants 1802319, 1807097, and 1952193.

\end{document}